\documentclass[a4paper,11pt]{amsart}

\usepackage[english]{my-shortcuts}
\usepackage{srcltx,algorithm,algorithmic}
\usepackage{a4wide}

\begin{document}

\title[Controllability of complex networks using perturbation theory]
{Controllability of complex networks using perturbation theory of extreme singular values}
\author{St\'ephane Chr\'etien} 
\author{S\'ebastien Darses}

\address{National Physical Laboratory\\ 
Hampton Road\\
Teddington, UK} 
\email{stephane.chretien@npl.co.uk}

\address{LATP, UMR 6632\\
Universit\'e Aix-Marseille, Technop\^ole Ch\^{a}teau-Gombert\\
39 rue Joliot Curie\\ 13453 Marseille Cedex 13, France \\
and \\
Laboratoire de Math\'ematiques, UMR 6623\\ 
Universit\'e de Franche-Comt\'e, 16 route de Gray,\\
25030 Besancon, France}
\email{sebastien.darses@univ-amu.fr}

\maketitle

\vspace{.5cm}

\begin{abstract}
	Pinning control on complex dynamical networks has emerged as a very important topic in 
	recent trends of control theory due to the extensive study of collective coupled behaviors
	and their role in physics, engineering and biology. In practice, real-world networks 
	consists of a large number of vertices and one may only be able to perform a control on a fraction of them only. Controllability of such systems has been addressed in \cite{PorfiriDiBernardo:Automatica08}, where it was reformulated as a global asymptotic stability problem. The goal of this short note is to refine the analysis proposed 
	in \cite{PorfiriDiBernardo:Automatica08} using recent results in singular value perturbation theory.
\end{abstract}


\section{Introduction}

In recent years, extensive efforts have been devoted to the control of complex dynamical networks. One major issue is that real world networks usually consist of a very large number of nodes and links which makes it impossible to apply control actions to all nodes. 

Pinning control is a new way to address this problem by placing local feedback injections on a small fraction of the nodes. 

Controllability of such systems has been addressed in \cite{PorfiriDiBernardo:Automatica08}, where it was reformulated as a global asymptotic stability problem. The goal of this short note is to refine the analysis proposed 
in \cite{PorfiriDiBernardo:Automatica08} using recent results in singular value perturbation theory.
\subsection{The model} 
One considers a set of 
$N$ $n$-dimensional oscillators governed by a system of nonlinear differential equations. 
Moreover, we assume that each oscillator is coupled with a restricted set of other oscillators. 
This coupling relationship can be efficiently described using a graph where the vertices 
are indexed by the oscillators and there is an edge between two oscillators if they are coupled. 
The overall dynamical system is given by the following set of differential equations
\bea \label{syst}
x_i^\prime(t) & = & f(x_i(t)) -\sigma B \sum_{j=1}^N l_{ij} x_j(t)+u_i(t), \ t\ge t_0,
\eea
$i=1,\ldots,N$, where $x_i(t) \in \R^n$ is the state of the $i^{th}$ oscillator, 
$\sigma>0$, $B\in \R^{n\times n}$, 
$f: \R \to \R$ describes the dynamics of each oscillator, $L=(l_{ij})_{i,j=1,\ldots,N}$ is the graph 
Laplacian of the underlying graph, and $u_i(t)$, $i=1,\ldots,N$ are the controls. 
For the system to be well defined, we have to specify some initial conditions $x_i(t_0)=x_{i0}$ for $i=1,\ldots,N$. 

\subsection{The control problem}
Assume that we have a reference trajectory $s(t)$, $t\ge t_0$ satisfying the differential equation 
\bean 
s^\prime(t) & = & f(s(t)). 
\eean 
Our goal is to control the system using a limited number of nodes. The selected nodes are called the "pinned nodes". 
For this purpose, we use a linear feedback law of the form 
\bean 
u_i(t) & = & p_i Ke_i(t), 
\eean 
where $e_i(t)=s(t)-x_i(t)$, $K$ is a feedback gain matrix, and where
\bean 
p_i & = & 
\begin{cases}
	1 \: \textrm{ if node $i$ is pinned} \\
	0 \: \textrm{ otherwise}.
\end{cases}
\eean 
Let $P$ denote the diagonal matrix with diagonal $(p_1,\ldots,p_N)$. 

\subsection{Controllability}
In \cite{PorfiriDiBernardo:Automatica08}, the authors propose a definition for (global pinning-) controllability (based on Lyapunov stability criteria):

\begin{defi}
	We say that the system (\ref{syst}) is controllable if the error dynamical system $e:=(e_i(t))_{1\le i\le N}$ is Lyapunov stable around the origin, i.e. there exists a positive definite function $V$ such that $\frac{d}{dt}V(e(t))< 0$ when $e(0)\neq 0$.
\end{defi}

The following result, \cite[Corollary 5]{PorfiriDiBernardo:Automatica08}, provides a sufficient condition for a system to be controllable:
\begin{prop}[\cite{PorfiriDiBernardo:Automatica08}]
	\label{cor5}
	Assume that $f$ is such that there exists a bounded matrix $F_{\xi,\wdt{\xi}}$, 
	whose coefficients depend on $\xi$ and $\tilde{\xi}$, which satisfies 
	\bea
	\label{Fxi} 
	F_{\xi,\wdt{\xi}}\left(\xi-\wdt{\xi}\right) & = & f(\xi)-f(\wdt{\xi}), \quad \xi,\wdt{\xi}\in\R^n.
	\eea  
	Let $Q\in \R^{n\times n}$ be a positive definite matrix such that 
	\bean 
	QK+K^tQ^t & = & \kappa \left(QB+B^tQ^t \right) \\
	\left(QB+B^tQ^t \right) & \succeq 0 
	\eean 
	and 
	\bea
	\label{tart}
	\frac12 \ \lb_{N} \left(\sigma L+\kappa P  \right)  \ \lb_{n}\left(QB+B^tQ^t \right) & > & \sup_{\xi,\wdt{\xi}}
	\|F_{\xi,\wdt{\xi}}\|\ \|Q\|.
	\eea 
	Then the system is controllable.
\end{prop}
Many systems of interest satisfy the constraint specified by (\ref{Fxi}); see \cite{JiangTangChen:CSF03}.
This proposition is very useful for node selection via the matrix $P$. Indeed, assume that $Q$ is selected, 
then one may try to maximise $\lb_{N} \left(\sigma L+\kappa P  \right)$ as a function of $P$, 
under the constraint that no more than $r$ nodes can be pinned. This is a combinatorial problem that 
can be relaxed using semi-definite programming or various heuristics \cite{GhoshBoyd:IEEECDC06}. 

\subsection{Goal of the paper}

Our goal in the present note is to propose an easy controllability condition refining \cite[Corollary 7]{PorfiriDiBernardo:Automatica08}, based on the algebraic connectivity of the graph, the number of pinned nodes, the coupling strengh and the feedback gain. Our approach is based on perturbation theory of the extreme singular values of a matrix after appending a column. The basic results of this theory are given in the appendix.

\section{Main result}

\subsection{Notations} The Kronecker symbol is denoted by $\delta_{i,j}$, i.e. $\delta_{i,j}=1$ if $i=j$ and
is equal to zero otherwise. We denote by $\|x\|_2$ the euclidian norm of a vector $x$ and by $\|A\|$ the associated operator norm (spectral norm) of a matrix $A$.

For any symmetric matrix $B\in \R^{d\times d}$ we will denote its eigenvalues by $\lb_1(B)
\ge \cdots \ge \lb_d(B)$. The largest eigenvalue will sometimes also be denoted by $\lambda_{\max}(B)$ and 
the smallest by $\lb_{\min}(B)$. The smallest nonzero eigenvalue of a positive semi-definite matrix $B$ will be denoted by $\lb_{\min>0}(B)$.

\subsection{A simple criterion for controllability}
Our main result is the following theorem. 
\begin{theo} \label{propkappa}
Let $Q\in \R^{n\times n}$ be a positive definite symmetric matrix that satisfies  
\bean 
QK+K^tQ^t & = & \kappa \left(QB+B^tQ^t \right) \\
\left(QB+B^tQ^t \right) & \succeq 0,
\eean
and assume that
\bea
\label{Fcond}
\|F_{\xi,\tilde{\xi}}\| & < & \frac{\sigma \lb_{\min>0}(L) \ \lb_{\min}\left(QB+B^tQ^t\right) }{2\ \|Q\|}.
\eea 
If $\kappa$ satisfies 
\bean 
\kappa & \ge & \frac{\sum_{i=1}^r {\rm deg}_i}{\sigma \lb_{\min>0}(L) - \frac{2\ \|F_{\xi,\tilde{\xi}}\| \ \|Q\|}{\lb_{\min}\left(QB+B^tQ^t\right)}}
+\sigma \lb_{\min>0}(L),
\eean 
then the system is controllable. 
\end{theo}
\begin{proof}
We follow the same steps as for the proof of Corollary 7 in \cite{PorfiriDiBernardo:Automatica08}. 
We assume without loss of generality that the first $r$ nodes are the pinned nodes. 
We may write $P$ as 
\bean 
P & = & \sum_{i=1}^r e_i e_i^t,  
\eean 
where $e_i$ is the $i^{th}$ member of the canonical basis of $\R^N$, i.e. $e_i(j)=\delta_{i,j}$. We will 
try to compare $\lb_{N} \left(\sigma L+\kappa P \right)$ with $\lb_{N} \left(\sigma L\right)$ and 
use Proposition \ref{cor5} to obtain a sufficient condition for controllability based on $L$, i.e. 
the topology of the network. For this purpose, let us recall that $L$ can be written as 
\bean  
L & = & \mathcal I\cdot \mathcal I^t,
\eean 
where $\mathcal I$ is the incidence matrix of any directed graph obtained from the system's graph by 
assigning an arbitrary sign to the edges \cite{BrouwerHaemers:Springer12}. 
Of course $L$ will not depend on the chosen assignment. Using this factorization of $L$, 
we obtain that 
\bean 
\sigma L + \kappa \sum_{i=1}^r  e_ie_i^t & = & 
\left[\sqrt{\kappa}  \ e_r,\ldots,\sqrt{\kappa} \ e_{1},\sqrt{\sigma} \mathcal I\right] \ 
\left[\sqrt{\kappa} \ e_r,\ldots,\sqrt{\kappa}  \ e_{1},\sqrt{\sigma} \mathcal I\right]^t.
\eean  
Moreover, $\lb_{\min>0}\left(\sigma L + \kappa P\right)$ can be expressed easily as the 
smallest nonzero eigenvalue of the $r^{th}$ term of a sequence of matrices with shape (\ref{A}) 
for which we can use Theorem \ref{smallestnonzero} iteratively. Indeed, we have 
\bean 
\lb_{\min>0}\left( \sigma L + \kappa  e_1 \right)
& = & 
\lb_{\min>0}\left( \left[\sqrt{\kappa} \ e_{1},\sqrt{\sigma} \mathcal I\right]^t
\left[\sqrt{\kappa} \ e_{1},\sqrt{\sigma} \mathcal I\right]\right).
\eean 
Let us denote by $x$ the vector $\sqrt{\kappa} \ e_1$ and by $X$ the matrix 
$[\sqrt{\sigma} \mathcal I]$. Then, we have that 
\bean 
\left[\sqrt{\kappa} \ e_{1},\sqrt{\sigma} \mathcal I\right]^t
\left[\sqrt{\kappa} \ e_{1},\sqrt{\sigma} \mathcal I\right] 
& = & 
\left[
\begin{array}{cc}
x^t x & x^t X \\
X^tx  & X^t X
\end{array}
\right].
\eean 
Therefore, Theorem \ref{smallestnonzero} gives 
\bean 
\lb_{\min>0} \left(\sigma L+\kappa e_1e_1^t \right) & \ge & \sigma \lb_{\min>0}(L) 
-\frac{{\rm deg}_1}{(\kappa-\sigma \lb_{\min>0}(L))},
\eean 
where ${\rm deg}_1$ is the degree of node number 1. 

Let us now consider 
$\lb_{\min>0} \left(\sigma L + \kappa \ e_1 +\delta_2 e_2 \right)$. We have that 
\bean 
\lb_{\min>0} \left(\sigma L + \kappa \ e_1 +\delta_2 e_2 \right)
& = & \lb_{\min>0} \left(
\left[\sqrt{\kappa} \ e_2,\sqrt{\kappa} \ e_{1},\sqrt{\sigma} \mathcal I\right]^t
\left[\sqrt{\kappa} \ e_2,\sqrt{\kappa} \ e_{1},\sqrt{\sigma} \mathcal I\right]\right).
\eean 
Let us denote by $x$ the vector $\sqrt{\kappa} \ e_2$ and by $X$ the matrix 
$[\sqrt{\kappa} \ e_1,\sqrt{\sigma} \mathcal I]$. Then, we have that 
\bean 
\left[\sqrt{\kappa} \ e_{2},\sqrt{\kappa} \ e_{1},\sqrt{\sigma} \mathcal I\right]^t
\left[\sqrt{\kappa} \ e_{2},\sqrt{\kappa} \ e_{1},\sqrt{\sigma} \mathcal I\right] 
& = & 
\left[
\begin{array}{cc}
x^t x & x^t X \\
X^tx  & X^t X
\end{array}
\right]
\eean
and using Theorem \ref{smallestnonzero} again, we obtain  
\bean 
\lb_{\min>0} \left(\sigma L+\kappa e_1e_1^t+\kappa e_2e_2^t \right) & \ge & \lb_{\min>0}(\sigma L+\kappa e_1e_1^t) 
-\frac{{\rm deg}_2}{(\kappa-\lb_{\min>0}(\sigma L+\kappa e_1e_1^t))}. 
\eean  
Since $\lb_{\min>0}(\sigma L+\kappa e_1e_1^t) \le \lb_{\min>0}(\sigma L)$, we thus obtain 
\bean 
\lb_{\min>0} \left(\sigma L+\kappa e_1e_1^t+\kappa e_2e_2^t \right) & \ge & \lb_{\min>0}(\sigma L+\kappa e_1e_1^t) 
-\frac{{\rm deg}_2}{(\kappa-\sigma \lb_{\min>0}(L))}. 
\eean  

We can repeat the same argument $r$ times and obtain 
\bea
\label{bnd>0} 
\lb_{\min>0} \left(\sigma L+\kappa P \right) & \ge & \sigma \lb_{\min>0}(L) -
\frac{\sum_{i=1}^r {\rm deg}_i}{\kappa-\sigma \lb_{\min>0}(L)}. 
\eea 

Finally, by Proposition \ref{cor5}, we know that the following constraint is sufficient for 
preserving controllability 
\bea
\label{condkappa}
\lb_{\min>0}\left(\sigma L+\kappa \sum_{i=1}^r e_ie_i^t\right) & \ge  & \frac{2\ \|F_{\xi,\tilde{\xi}}\| \ \|Q\|}{\lb_{\min}\left(QB+B^tQ^t\right)}.
\eea
By (\ref{bnd>0}), it is sufficient to garantee the controllability of our system to impose 
\bean 
\sigma \lb_{\min>0}(L) -
\frac{\sum_{i=1}^r {\rm deg}_i}{\kappa-\sigma \lb_{\min>0}(L)} & \ge & \frac{2\ \|F_{\xi,\tilde{\xi}}\| \ \|Q\|}{\lb_{\min}\left(QB+B^tQ^t\right)}.
\eean 
Then, combining (\ref{condkappa})  with (\ref{Fcond}) implies that 
\bean 
\kappa & \ge & \frac{\sum_{i=1}^r {\rm deg}_i}{\sigma \lb_{\min>0}(L) - \frac{2\ \|F_{\xi,\tilde{\xi}}\| \ \|Q\|}{\lb_{\min}\left(QB+B^tQ^t\right)}}
+\sigma \lb_{\min>0}(L)
\eean 
is a sufficient condition for controllability.
\end{proof}

\appendix

\section{Perturbation theory of extreme singular values after appending a column}

\subsection{Framework}
Let $d$ be an integer. Let $X\in\R^{d\times n}$ be a $d\times n$-matrix and let $x\in \R^{d}$ be column vector.  We denote by a subscript $^t$ the transpose of vectors and matrices. 
There exist at least two ways to study the singular values of the matrix $(x,X)$ obtained by appending the column vector $x$ to the matrix $X$:
\begin{enumerate}
	\item[{\bf (A1)}] Consider  the matrix 
	\bea \label{add}
	A & = & 
	\left[
	\begin{array}{c}
		x^t \\
		X^t
	\end{array}
	\right]
	\left[
	\begin{array}{cc}
		x & X
	\end{array}
	\right]
	=
	\left[
	\begin{array}{cc}
		x^tx & x^tX \\
		X^tx & X^tX
	\end{array}
	\right];
	\eea
	\item[{\bf (A2)}] Consider  the matrix 
	\bean
	\wdt{A} & = & 
	\left[
	\begin{array}{cc}
		x & X
	\end{array}
	\right]\left[ 
	\begin{array}{c}
		x^t \\
		X^t
	\end{array}
	\right]
	= XX^t+xx^t.
	\eean 
\end{enumerate}

On one hand, one may study in {\bf (A1)} the eigenvalues of the $(n+1)\times (n+1)$ hermitian matrix $A$, i.e. the matrix $X^tX$ augmented with an arrow matrix.

On the other hand, one will deal in {\bf (A2)}  with the eigenvalues of the $d\times d$ hermitian matrix $\wdt A$, which may be seen as a rank-one perturbation of $XX^t$. The matrices $A$ and $\wdt A$ have the same non-zeros eigenvalues, and in particular 
$\lb_{\max}(A) = \lb_{\max}(\wdt A)$.
Moreover, the singular values of 
the matrix $(x,X)$ are the square-root of the eigenvalues of the matrix $A$.

Equivalently, the problem of a rank-one perturbation can be rephrased as 
the one of controlling the perturbation of the singular values of a matrix after appending a column. 

\subsection{A theorem of Li and Li} 
In this paper, we use a slightly more general framework than  {\bf (A1)}, that is the case of a matrix
\bea
\label{A}
A & = & 
\left[
\begin{array}{cc}
	c & a^t \\
	a & M
\end{array}
\right],
\eea
where $a\in \R^{d}$, $c\in \R$ and $M\in \R^{d\times d}$ is a symmetric matrix.

The following theorem provides sharp upper bounds for $\lb_{\max}(A)$, 
and lower bounds on $\lb_{\min}(A)$, depending on various information on the sub-matrix $M$ of $A$. As discussed 
above, this problem has close relationships with our problem of appending a column to a given rectangular matrix,
because $\lb_1(\wdt{A})=\lb_1(A)$. 

\begin{theo}[Li-Li's inequality and a lower bound]
	\label{main}
	Let $d$ be a positive integer and let $M\in \C^{d\times d}$ be an Hermitian matrix, whose eigenvalues are $\lambda_1\ge \cdots \ge \lambda_{d}$ with corresponding eigenvectors $(V_1,\cdots,V_d)$. Set $c\in \R$, $a\in \C^{d}$.
	Let $A$ be given by (\ref{A}). Therefore:
	\beq \label{lbmax}
	\frac{2\la a,V_1\ra^2}{\eta_1+\sqrt{\eta_1^2+4\la a,V_1\ra^2}} \le \lb_{1}(A)-\max(c,\lb_1) \le  \frac{2\|a\|^2}{\eta_1+\sqrt{\eta_1^2+4\|a\|^2}},
	\eeq
	with
	\bean
	\eta_1 & = & |c-\lb_1|.
	\eean
\end{theo}

\subsection{Perturbation of the smallest nonzero eigenvalue}

The same technics used to prove Theorem \ref{main} also give lower bounds for the smallest nonzero eigenvalue, which are also direct consequences of Li-Li's inequality. For more details, we refer the reader to \cite{ChretienDarses:ExtremePertSurvey}.

\begin{theo}
	\label{smallestnonzero}
	Let $d$ be a positive integer and let $M\in \C^{d\times d}$ be a positive semi-definite Hermitian matrix, whose eigenvalues are $\lambda_1\ge \cdots \ge \lambda_{d}$ with corresponding eigenvectors $(V_1,\cdots,V_d)$. Set $c\in \R$, $a\in \C^{d}$.
	Let $A$ be given by (\ref{A}). Assume that $M$ has rank $r\le d$. Therefore:
	\beq \label{lbmin}
	\lb_{r+1}(A) \ge \min(c,\lb_r)- \frac{2\|a\|^2}{\eta_r+\sqrt{\eta_r^2+4\|a\|^2}},
	\eeq
	with
	\bean
	\eta_r & = & |c-\lb_r|.
	\eean
\end{theo}

In particular, the following perturbation bounds of Weyl and Mathias hold:

\begin{coro}
	\bea
	\lb_{r+1}(A) & \ge &\min(c,\lb_{r}) - \|a\|_2 \label{ineg1min}\\
	\lb_{r+1}(A) & \ge &\min(c,\lb_{r}) - \frac{\|a\|_2^2}{|c-\lb_{r}|}.   \label{ineg2min}
	\eea
\end{coro}

{\bf Acknowledgement --} We thank Maurizio Porfiri for pointing out a missing assumption in 
Theorem \ref{propkappa}.



\bibliographystyle{amsplain}

\end{document}